\documentclass[12pt,titlepage]{article}
\usepackage{graphicx}

\begin{document}

\begin{center}
{\Large A new approach for analyzing panel AR(1) series
with application  to the unit root test}
\end{center}

\begin{center}
 Yu-Pin Hu
\footnote {Corresponding author. Professor of Department of International
Business Studies, National Chi Nan University, Taiwan. Address: Department of
International Business Studies, National Chi Nan University, No. 1 University
Road, Puli, Nantou 545, Taiwan.  Email: huyp@ncnu.edu.tw}
  and J. T. Gene Hwang
 \footnote { Chair professor of  Department of Mathematics, National Chung Cheng University, Taiwan.
 Professor Emeritus of Department of Mathematics, Cornell University, USA.
Email: hwang@math.cornell.edu
}
 \end{center}
\begin{center}
 {\footnotesize May 21, 2015 }\\
  \vspace{0.1in}
SUMMARY
\end{center}
This paper derives several novel tests to improve on the t-test for
testing AR(1) coefficients of panel time series, i.e., of multiple
time series, when each has a small number of observations. These
tests can determine the acceptance or the rejection of each
hypothesis individually while controlling the average type one
error. Strikingly, the testing statistics derived by the empirical
Bayes approach can be approximated by a simple form similar to the
t-statistic; the only difference is that the means and the variances
are estimated by shrinkage estimators. Simulations demonstrate that
the proposed tests have higher average power than the t-test in all
settings we examine including those when the priors are
miss-specified and the cross section series are dependent.
 \\
{}\\
{JEL Classification: C12; C32.}\\
 {\bf Keywords}:{ Empirical Bayes; Multiple tests; Panel time series; Random effect model;
  Shrinkage estimator. }

\newpage

\section{Introduction}
Testing a unit root hypothesis is a very important subject since, for example, it can be applied to test
 the Purchasing Power Parity (PPP) theory in Economics.
The topic attracts much attention in research also because the traditional unit
root tests can have low power in some circumstances. One such  circumstance
relates to the panel time series which consists of $N$ series to be
simultaneously tested, each having $T$ observations, where $N$ is large or
moderate and $T$ is small.  This is the scenario we consider in this paper.
Although we focus on the Economics settings when discussing applications, our
solutions are applicable to other applications involving the ``small $n$ and
large $p$" problems in Statistics.

In  this paper, we consider a new approach to derive several novel tests  that
 improve, in power, on the traditional t-test for testing multiple unit root
null hypotheses and white noise null hypotheses. This new approach is based on
the optimal   multiple test criterion (Liu, 2006, Storey, 2007, Storey et al.,
2007, Hwang and Liu, 2010, and Noma and Matsui, 2012 and 2013), which are
motivated by biological microarray data analyses.  This approach
 could determine acceptance or rejection of each
 hypothesis  individually while controlling the average type one error of the multiple tests.
Traditionally,  the panel unit root approach tests against  the null hypothesis
that all series follow a unit root model  (Levin et al., 1992, Baltagi and Kao,
2000, Bai and Ng, 2004, 2010, Pesaran, 2007, Pesaran, et al., 2013, and etc.).
Hence either all series are declared a unit root model (i.e. non-stationary
series) or some declared stationary series without identifying which. In
contrast, the tests proposed in this paper could determine the stationarity of
each series and in the meantime control the average type one error. This seems
more desirable.

The multiple test criterion considered by Liu (2006),
 Storey (2007), Hwang and Liu (2010), and Noma and Matsui (2012, 2013),
 are to maximize the average power while controlling the average type one error.
Interestingly, such a criterion is equivalent to other optimality criteria
based on controlling the false discovery rate (FDR). See Storey (2007) and
Hwang and Liu (2010). Their approaches and the approach in this paper all use
the Neymann-Pearson fundamental lemma to derive optimum procedures.

While all aim at controlling the average Frequentist type one errors, the
difference between the approaches of Storey (2007), and the group of
researchers,  Hwang and Liu (2010) and Noma and Matusi (2012, 2013) is that
Storey's approaches aims to maximize the Frequentist average powers whereas the
others aim to maximize the Bayesian average powers.  The advantage of the
approaches of Hwang and Liu (2010) and Noma and Matusi (2012, 2013) over
Storey's approach is that the former are much faster in computation and also,
as shown in Hwang and Liu (2010), provide higher average test power. While the
approaches of Hwang and Liu (2010) and Noma and Matusi (2012, 2013) are
similar, the statistics proposed by Hwang and Liu have further approximated
formulae in simpler forms, which can be easily calculated without evaluating
integral, unlike Norma and Matsui's  approach which requires  evaluating $3N$
integrals.

In this paper, we tackle the difficult problem of testing
coefficients of time series models. We follow the approach of Hwang
and Liu (2010) which constructs the MAP test, i.e., the test that
maximizes the Bayesian expected average power with respect to a
prior distribution  while controlling the Frequentist average type
one error. The general theory developed in Section 3 shows that the
MAP test is an approximation of Story's test.  To derive the
statistics   for testing the one-sided and two-sided hypotheses of
the coefficients of panel AR(1) models, we assume a class of priors
on the means, a class of priors on the variances, or on both
resulting in a MAP statistic shrinking the means, the variances or
both respectively. These statistics are further approximated,
leading to the proposed statistics.  Strikingly, in all situations,
the proposed statistics basically take a simple form similar to the
t-statistic; the only difference is that the means and the variances
are estimated by shrinkage estimators. Previously, such a result was
available in Hwang and Liu (2010), Cui er al. (2005), and Smyth
(2004) only for the usual ANOVA models and only for the procedure
shrinking the variances. For the procedure shrinking the means and
the variances, the tests of Hwang and Liu (2010) were not put in the
form of the t-statistic.

Our proposed shrinkage t-tests are  shown to have higher average powers than
the traditional t-test because the tests ``borrow the strength'' from all
series. The tests implicitly determine how similar the parameters of the series
are. The more similar
 the series are, the  more extensively  the data from other cross sections are used to estimate the
 parameters of the individual series.  Consequently, the improvements are larger.

Note that the testing statistics developed in this paper aim to satisfy the
Frequentist criteria of controlling the average
 type one error, even though
 we use a Bayesian approach to construct the proposed statistics.
To be more realistic, we consider  not only a prior but a class of priors
indexed by some hyper-parameters. We use the data to estimate the
hyper-parameters; hence the procedure  is called empirical Bayes, which  is
equivalent to the Frequentist approach based on a random effect model.
 Hence our results are quite
different from the  Bayesian unit root tests proposed in Uhlig (1994) and
Philips and Xiao (1998).

This paper uses a bootstrap method  to obtain the critical value to control the
average type one error of the proposed tests. Simulations in Section 7 show
that bootstrap works well for all our settings. Note that our problem is
different from the unit root bootstrap tests proposed by  Ferretti and Romo
(1996) and Park (2003), which aim at testing a single hypothesis with a large
number, $T$, of observations.  Simulation results also show that the proposed
tests have either higher or similar average power when compared with the
t-tests in all the cases we considered.  Specifically, when $N=80$ and $T=10$,
 the proposed tests increase the average power of  t-test
by $70 \%$ and $25 \%$, respectively,  for testing the white noise null
hypotheses and
  the unit root null hypotheses.
 We also demonstrate similar improvement  when the model is
misspecified and when the cross section series are dependent. In this paper,
although we only work on AR(1) models, we anticipate that these results can be
generalized to more complex time series models.

Our proposed tests are fast in computation. Given a data set with $N=1000$ and
$T=10$, it takes  about 10 seconds to compute our proposed tests
($F_{\mbox{ss}}$ and $RF_{\mbox{ss}}$) for all 1000 hypotheses, using a Laptop
and the GAUSS 9.0 program.

 The rest of this paper is organized as follows. In Section 2,
 we present the model considered in this research and give
a review of the optimal discovery procedure (Storey, 2007) and the maximizing
average power test (MAP) (Hwang and Liu, 2010). In Section 3, we develop a
theory that links the two approaches. In Sections 4 and 5, we derive the MAP
test under various prior assumptions for the  two-sided and the one-sided
hypotheses. The proposed empirical Bayes tests are constructed in Section 6,
where the issues of estimating the hyper-parameters and controlling the average
type one error by bootstrap are discussed. In Section 7, we present the
simulation results. Section 8 gives the concluding remarks.

\section{The Model and Reviews of  the ODP and  MAP Tests}

Suppose the $N$-dimensional AR(1) processes are generated by
\begin{equation}  \label{model-1}
y_{j,t}=\phi_j y_{j,t-1}+e_{j,t} \quad \mbox{for $1 \leq j \leq N$ and $1 \leq t \leq T$},
\end{equation}
where for section $j$, $e_{j,t}$ is an i.i.d. normal random variable
with zero mean and variance $\sigma_j^2$. Note that we allow the
dependence of the cross section series in model (\ref{model-1}).
Except in Sections 6 and 7, we derive the tests without assuming
that the cross section series are independent throughout the paper.

The observation of the $j$-th section is ${\bf y}_j=(y_{j,1},\cdots,y_{j,T}
)'$, and the probability density function (pdf) of ${\bf y}_j$ given $y_{j,1}$
is
\begin{equation}   \label{pdf-1}
 (\frac{1}{\sqrt{2\pi}\sigma_j})^{T-1}
e^{\frac{-1}{2} \sum^T_{t=2} (y_{j,t}-\phi_j y_{j,t-1})^2/\sigma^2_j }.
\end{equation}
It is easy to see that $\sum^T_{t=2} (y_{j,t}-\phi_j y_{j,t-1})^2 =
\sum^T_{t=2} (y_{j,t}-{\hat \phi}_j y_{j,t-1})^2 + ({\hat
\phi}_j-\phi_j)^2\sum^T_{t=2} y_{j,t-1}^2 $, where
\begin{equation} \label{ols}
{\hat \phi}_j =\frac{\sum^{T}_{t=2} y_{j,t} y_{j,t-1}}{\sum^{T}_{t=2} y^2_{j,t-1}}
\end{equation}
is the least squares estimator by regressing ${y_{j,t}}'s$ against
${y_{j,t-1}}'s$. Using the notation
\begin{equation} \label{sigma}
{\hat \sigma_j^2}=(1/(T-1)) \sum^{T}_{t=2}(y_{j,t}-{\hat \phi}_j y_{j,t-1})^2 \quad
\mbox{and} \quad
\quad S_{j}= \sum^{T}_{t=2} y^2_{j,t-1},
\end{equation}
the pdf of ${\bf y}_j$ can be written as
\begin{equation} \label{pdf-2}
 (\frac{1}{\sqrt{2\pi}\sigma_j})^{T-1}
{e^{-\frac{1}{2} (T-1) {\hat \sigma}^2_j/\sigma^2_j }}
\cdot
{e^ {-\frac{1}{2} ({\hat \phi}_j-\phi_j)^2S_j/\sigma^2_j }}.
\end{equation}
It follows that $({\hat \phi}_j,{\hat \sigma}^2_j, S_j)$ is a sufficient statistic. Also
$({\hat \phi}_j,{\hat \sigma}^2_j)$ is the maximum likelihood estimator of
$({\phi}_j,{\sigma}^2_j)$.

In this paper, we consider testing simultaneously  the $N$ hypotheses
\begin{equation}  \label{hyp}
H^j_0: \quad \phi_j =\phi_0 \quad \mbox{vs.} \quad H^j_1: \quad \phi_j \in D,
\end{equation}
where $\phi_0$ is a fixed number, $j=1,\cdots,N$, and $D$ denotes a set of
${\phi_j}'s$. When  $D$ consists of a single point, (\ref{hyp}) corresponds to
a simple test. When
 $D=\{\phi| \phi \neq \phi_0 \}$, (\ref{hyp}) corresponds to
 a two-sided test and when
 $D=\{\phi| \phi \leq \phi_0 \}$, (\ref{hyp}) corresponds to
 a   one-sided test.

The optimal discovery procedure (ODP) of Storey (2007) aims at constructing a
rejection region that maximizes the expected number of the true positive
results (ETP) while controlling  the expected number of false positive results
(EFP). The discussion below is applicable to the general situation, where ${\bf
y}_j$, a $T-$dimensional random vector, has  pdf $f({\bf
y}_j|\phi_j,\sigma_j)$. The procedures we consider are to
\begin{equation}\label{rejectc}
\mbox{reject $\quad$ $H^j_0$ $\quad$ if $\quad$ ${\bf y}_j \in C$,} \end{equation}
 where $C$, independent of $j$, is a set of $T$-dimensional
vectors.  Storey (2007) argues that procedures using $C$ depending on $j$ have
no advantage in average power. Using $I(\cdot)$ to denote  an indicator
function, we have {\small
\begin{equation}  \label{etp-1}
\mbox{ ETP}= E \{\sum^N_{j=1} I(\phi_j \in D \quad \mbox{\small and} \quad {\bf y}_j \in C) \}
=\sum_{\{j| \phi_j \in D \}} P_{\phi_j,\sigma_j} ( {\bf y}_j \in C)
=\int_{y \in C} \sum_{\{j| \phi_j \in D\}}  f({\bf y}|\phi_j,\sigma_j) d {\bf y}.
\end{equation}
}
Similarly,
\begin{equation}  \label{efp-1}
\mbox{ EFP}= \sum_{\{j|\phi_j =\phi_0\}} P_{\phi_j,\sigma_j} ( {\bf y}_j \in C) =
 \int_{y \in C} \sum_{\{j|\phi_j = \phi_0\}}  f({\bf y}|\phi_j,\sigma_j) d {\bf y}.
\end{equation}

Let $N_1$ and $N-N_1$ be the number of series satisfying the  alternative and
the null hypotheses, respectively. Note that $\mbox{ETP}/N_1$ and
$\mbox{EFP}/(N-N_1)$ are therefore the average power and the average type one
error, which are quantities of great concern to a Frequentist.

Applying Neymann-Pearson lemma to (\ref{etp-1}) and  (\ref{efp-1}) leads to the
 optimal test which maximizes ETP
(or equivalently $\mbox{ETP}/N_1$) while controlling $\mbox{EFP}/(N-N_1)$ at
level $\alpha$. Assume that the pdf of ${\bf y}_j$ is $f({\bf
y}_j|{\phi_j,\sigma_j})$, where ${\phi_j}'s$ are the key parameters whereas
${\sigma_j}'s$ are the nuisance parameters, which can be interpreted as
variances. See, for example, the pdf in (\ref{pdf-1}). The rejection region is
then ${\bf y}_j \in C$, where
\begin{equation}  \label{rej-1}
C=\{{\bf y} |\quad \frac{ \sum_{ \{j|\phi_j \in D \}}  f({\bf y} |{\phi_j,\sigma_j})
 }{  \sum_{\{j|\phi_j =\phi_0\}}  f({\bf y}|{\phi_j,\sigma_j })
} > \mbox{crit}  \},
\end{equation}
where crit is a cutoff point chosen so that it has average type one error equal to $\alpha$.

There is, however, a problem with the ``test" in  (\ref{rej-1}). In order to
apply  (\ref{rej-1}), one needs to know which hypothesis is true and which is
false, the very information that one is trying to determine. When $\sigma_j$
are all equal (to $\sigma$), there is however a possible way out, as described
in Storey (2007). Write the inequality  in (\ref{rej-1}) as $A/D> \mbox{crit}$
where $A$ and $D$ represent the numerator and the denominator of sums of
probability density functions. Notice that the inequality in (\ref{rej-1}) is
equivalent to $(A+D)/D >\mbox{crit}+1$. Also $D$ equals $(N-N_1) \cdot f({\bf
y}|{\phi_0,\sigma})$. Putting all these together and omitting some constants,
the statistic is equivalent to
 $\sum^N_{j=1} f({\bf y}|{\phi_j,\sigma})/f({\bf y}|{\phi_0,\sigma})$, which can be calculated without knowing
   which hypothesis is true or false.

However, when ${\sigma_j}'s$ are different, it is much harder to approximate the statistic
in  (\ref{rej-1}).
Storey, et al. (2007) did have a successful attempt, where  ${\sigma_j}$ is replaced by an estimator
 based on the $j$th population. Storey's procedure, however, is computationally very intensive.
  It requires calculating
$N$ times the statistic, which involves summation of $N$ terms. When $N$ is
large, it is overwhelming.

The approach of Hwang and Liu (2010) is more parametric because of postulating
classes of prior distributions on both $\phi_j$ and $\sigma_j$. They construct
their MAP (acronym of maximum average power) test to maximize the Bayesian
expected value of ETP, which is the average power with respect to the prior
distribution. Using some intuitive approximation in the empirical Bayes
fashion, they were able to derive some statistics which can be calculated
instantaneously. In particular, their approach leads to a statistic, called
$F_{\mbox{ss}}$, not only borrowing the strength from all populations to
estimate $\phi_j$ (which Storey's procedure does), but also to estimate
$\sigma_j$ (which Storey's procedure does not). Consequently, it is to be
expected that $F_{\mbox{ss}}$ test has higher average power, which was
numerically demonstrated to be so.

In this paper, the approach of Hwang and Liu (2010) is applied to the time
series models which are much more difficult to construct statistical tests than
their ANOVA models. To the best of our knowledge,  this paper is the first to
present shrinkage multiple tests for the time series model. The empirical Bayes
approach  is equivalent to the random effect model approach, because  the
parameters are assumed to be random in either case.

\section{The Main Theorems}
In this section, we shall provide a general theory that shows that the ODP
approach of Storey (2007) is asymptotically equivalent to the MAP test of Hwang
and Liu (2010). We then apply the theory to the AR(1) model in the following
sections. In the case when a theorem needs a proof, the proof is provided in
the Appendix. In this section, we consider testing the hypothesis (\ref{hyp})
by assuming that each ${\bf y}_j$ has the pdf $f({\bf y}_j|\phi_j,\sigma_j)$.
Under the assumption, we consider three cases of priors:
\begin{description}
\item{{\bf Case 1:} ${\phi}_j$ is fixed and  ${\sigma}_j's$  are i.i.d. each having the distribution $\pi_1(\sigma)$.}
\item{{\bf Case 2:} $\sigma_j's$ are fixed and  $\phi_j$ are i.i.d., each having the distribution $\pi_2(\phi)$.}
\item{{\bf Case 3:} $(\phi_j,\sigma_j)$ are i.i.d., each having the distribution $\pi(\phi,\sigma)$.}
\end{description}

Considerations of these three cases shall, in order, lead to  $F_{\mbox{sv}}$,
$F_{\mbox{sm}}$ and $F_{\mbox{ss}}$ tests for the two-sided hypothesis and to
 $RF_{\mbox{sv}}$,  $RF_{\mbox{sm}}$ and $RF_{\mbox{ss}}$ for the one-sided hypothesis.
The subscripts sv, sm, and ss, represent tests that shrink the variances,
shrink the means and shrink both the means and variances, respectively.

We first write an asymptotic formula for a rejection region $C$ for Case 1,
which rejects $H^j_0$ if and only if ${\bf y}_j \in C$.
\begin{description}
\item{\bf Theorem 1 (Case 1).} Under the assumption of Case 1, as $N_1$ and $N-N_1$ go to infinity,
\begin{equation} \label{3-1}
\frac{\mbox{ETP}}{N_1}-\mbox{BETP} \longrightarrow 0 \quad \mbox{in probability}
\end{equation}
and
\begin{equation} \label{3-2}
\frac{\mbox{EFP}}{N-N_1}-\mbox{BEFP} \longrightarrow 0 \quad \mbox{in probability},
\end{equation}
where $$ \mbox{BETP}= \frac{1}{N_1} \sum_{\{j|\phi_j \in D\}} \int
P_{\phi_j,\sigma}({\bf y} \in C) d \pi_1 (\sigma),$$ and
 $$
\mbox{BEFP}= \int P_{\phi_0,\sigma}({\bf y} \in C)
 d \pi_1(\sigma).
$$
\end{description}

The ``B" in the notation of BETP and BEFP stands for ``Bayes". Actually, both
quantities are also the Bayesian expected values of  the Frequentist's
quantities,  ${\mbox{ETP}}/N_1$ and ${\mbox{EFP}}/(N-N_1)$. The MAP test of
Hwang and Liu (2010) is defined as the test that maximizes the Bayesian
expectation
 of the average power, BETP, among all tests such that
\begin{equation} \label{3-3}
\mbox{BEFP} \leq \alpha.
\end{equation}
Hence the theorem shows that the ODP approach of Storey (2007) is
asymptotically equivalent to the MAP test of Hwang and Liu (2010).

 The proof of the above theorem using the law of large numbers is based on
the assumption that $\sigma_j's$ are independent. However, even if $\sigma_j's$
are correlated, it is possible to write weaker conditions so that the law of
large numbers applies and hence the theorem could be established under weaker
assumptions.

Next, by interchanging the order of integration, we can write BETP and BEFP as
$$ \mbox{BETP} =
\int_{{\bf y} \in C}  \frac{1}{N_1} \sum_{\{j|\phi_j \in D\}} \int f({\bf y}| {\phi_j, \sigma})  d\pi_1 (\sigma) d{\bf y}
$$
and
$$
\mbox{BEFP}= \int_{{\bf y} \in C}  \int f({\bf y}| {\phi_0, \sigma})  d\pi_1 (\sigma) d{\bf y}.
$$
Therefore  the Neymann-Pearson fundamental lemma implies the following theorem.
\begin{description}
\item{\bf Theorem 2 (Case 1).} Among all the procedures in (\ref{rejectc}),
 the MAP test consists of ${\bf y}$ such that
\begin{equation} \label{3-4}
\mbox{PT}(\phi_1,\cdots,\phi_N) \equiv \frac{
 \sum_{\{j|\phi_i \in D\}} \int f({\bf y}|\phi_j,\sigma) d \pi_1(\sigma)}
{ \int f({\bf y}|\phi_0,\sigma) d \pi_1(\sigma)} >\mbox{crit},
\end{equation}
where crit is  chosen so that equality in (\ref{3-3}) is attained for the test  (\ref{3-4}).
\end{description}

There is, however, a problem with PT in  (\ref{3-4}), which stands for ``pseudo
test". Namely, it still depends on the unknown parameters $\phi_j's$ and is not
really an applicable test. Following the principle of likelihood ratio test, we
can use the statistic
\begin{equation} \label{3-5}
\mbox{sup}_{\phi_1,\cdots,\phi_N \in D}\mbox{PT}(\phi_1,\cdots,\phi_N),
\end{equation}
which leads to a real statistical test. Note we can view  (\ref{3-5}) as an
approximate MAP test, since the $\phi_j's$ are replaced by the maximum
likelihood estimators ${\hat \phi}_j's$. The approximation is one of the best
imaginable approximations, even when the sample sizes are small.

We now state the theorem which, under a condition, gives us an explicit formula
for (\ref{3-5}).
\begin{description}
\item{\bf Theorem 3 (Case 1).}
Assume that the maximization of $f({\bf y}|\phi,\sigma)$ with respect to $\phi \in D$ is attained
when $\phi={\hat \phi}_M ({\bf y}) \in D$ and ${\hat \phi}_M$ does not depend on $\sigma_j$. Then
(\ref{3-5}) equals
\begin{equation} \label{3-6}
\frac{\int f({\bf y}| {\hat \phi}_M, \sigma) d \pi_1(\sigma)}
{\int f({\bf y}| \phi_0, \sigma) d \pi_1(\sigma)}.
\end{equation}
\end{description}

For Case 2, similar to Theorems 1 and 2, we could  establish results
 which are stated below while omitting the proof.

\begin{description}
\item{\bf Theorem 4 (Case 2).} Under the assumption of Case 2,
the statement  in Theorem 1 holds with
  $$
\mbox{BETP}= \frac{1}{N_1} \sum_{\footnotesize{\{\mbox{$j| H^j_0$ is false} \}}}
 \int  P_{\phi,\sigma_j} ({\bf y} \in C)
 d \pi_2(\phi)
$$
and
 $$
\mbox{BEFP}= \frac{1}{N-N_1} \sum_{\footnotesize{\{\mbox{$j| H^j_0$ is true}\}}}
 P_{\phi_0,\sigma_j} ({\bf y} \in C).
$$
Also the MAP test for this case rejects $H^j_0$  if and only if ${\bf y}_j \in C$ and
\begin{equation} \label{3-8}
C=\{{\bf y}| \frac{  \sum_{\footnotesize{\{\mbox{$j| H^j_0$ is false} \}}}
 \int  f ({\bf y}|{\phi,\sigma_j}) d \pi_2(\phi)}{
   \sum_{\footnotesize{\{\mbox{$j| H^j_0$ is true}\}}}
 f ({\bf y}|{\phi_0,\sigma_j}) }
 > \mbox{crit}
 \},
\end{equation}
where crit is chosen so that the equality in  (\ref{3-3}) holds.
\end{description}

The above region is not a usable rejection region, since it depends on unknown
$\sigma_j's$. To derive a useful version, we could use the likelihood ratio
principle by taking the sup of the numerator and denominator of the ratio in
(\ref{3-8}). It is easy to see that the resultant ratio is proportional to
\begin{equation} \label{3-9}
 \frac{  \int  f ({\bf y}|{\phi,{\hat \sigma}_M}) d \pi_2 (\phi)}{
  f ({\bf y}|{\phi_0,{\hat \sigma}_0})},
\end{equation}
where ${\hat \sigma}_M$ maximizes  $\int  f ({\bf y}|{\phi,{\sigma}}) d \pi_2
(\phi)$ and ${\hat \sigma}_0$ maximizes  $\int  f ({\bf y}|{\phi_0,{\sigma}}) d
\pi_2 (\phi)$, since
 ${\hat \sigma}_M$ and  ${\hat \sigma}_0$ do not depend on $j$.
 Now (\ref{3-9}) is a real statistic, since it does not depend on $\sigma_j's$.
    We note here that the process of turning (\ref{3-8}) into a real statistic can also be carried out
    with estimators which are different  from the maximizers. Later on, it turns out that
 the close form  of  ${\hat \sigma}_M$ can not be easily derived for the AR(1) model and so
we use a different estimator ${\hat \sigma}_{\star}$ to substitute for
$\sigma_j$ in (\ref{3-8}). This leads to (\ref{3-9}) with  ${\hat \sigma}_M$
being replaced by  ${\hat \sigma}_{\star}$ i.e.
\begin{equation} \label{3-10}
 \frac{  \int  f ({\bf y}|{\phi,{\hat \sigma}_{\star}}) d \pi_2 (\phi)}{
  f ({\bf y}|{\phi_0,{\hat \sigma}_0})}.
\end{equation}

We finally came to the easiest case, Case 3. We state the following theorem and omit
the proof, which is similarly to that of Theorem 1.

\begin{description}
\item{\bf Theorem 5 (Case 3).} Under the assumption of Case 3,
the statement  in Theorem 1 holds with
$$
\mbox{BETP}=
 \int  P_{\phi,\sigma} ({\bf y} \in C)
 d \pi(\phi,\sigma)
 =  \int_{{\bf y} \in C} \int  f({\bf y}|\phi,\sigma)
 d \pi(\phi,\sigma) d {\bf y}
$$
and
$$
\mbox{BEFP}=  \int
 P_{\phi_0,\sigma} ({\bf y} \in C)  d \pi_1(\sigma)
  =  \int_{{\bf y} \in C} \int  f({\bf y}|\phi_0,\sigma)
 d \pi_1(\sigma) d {\bf y}.
$$
  Consequently,
the MAP test rejects $H^j_0$ if ${\bf y}_j \in C$ and
 \begin{equation}  \label{map}
C = \{{\bf y} |\quad \frac{
 \int_{ \phi \in D}  f ({\bf y}| \sigma \mbox{,  } \phi)   d \pi(\phi,\sigma)}
 {  \int f ({\bf y} | \sigma \mbox{,  }\phi_0)    d \pi_1(\sigma)  } > \mbox{crit} \},
 \end{equation}
where crit is chosen so that the equality in  (\ref{3-3}) holds.
\end{description}

The likelihood ratio principle is not used here in deriving Theorem 5.

\section{The Two-sided Test}
\subsection{The t-test}
To begin, we consider the two-sided tests
\begin{equation}
H^j_0: \quad \phi_j =\phi_0 \quad \mbox{vs.} \quad H^j_1: \quad \phi_j \neq \phi_0, \quad
\mbox{where $j=1,\cdots,N$}.
\end{equation}
The well-known  t-test is to reject $H^j_0$ if $t^2_j$ is larger than a
critical value, where
\begin{equation} \label{t-test}
t_j= ({\hat \phi}_j -\phi_0) (\frac{S_{j}} {\hat \sigma_j^2})^{1/2}.
\end{equation}

This test is asymptotically optimal in power if we consider each hypothesis separately. However,
tests with larger average power can be constructed as outlined in Section 3.
 Similar to Hwang and Liu (2010), we construct $F_{\mbox {sv}}$, $F_{\mbox{sm}}$, and $F_{\mbox{ss}}$,
which shrink the variances (sv), the means (sm) and  both the variances and the
means (ss), respectively.

\subsection{The test shrinking the variances: $F_{\mbox{sv}}$ }
To shrink only the variances and not the means, we shall consider  Case 1 in
Section 3. The pdf  $f({\bf y}_j|\phi_j,\sigma_j)$ is given in  (\ref{pdf-2})
and $D=\{\phi: \phi \neq \phi_0 \}$.  From  (\ref{pdf-2}), note that the
condition in Theorem 3 is satisfied with ${\hat \phi}_M={\hat \phi}$, where
${\hat \phi}$ is defined in (\ref{ols}) unless ${\hat \phi}_M={\phi}_0$. The
latter situation occurs with zero probability and hence can be ignored. It may
be easier for the future user to have formulas with ${\bf y}$ replaced with
${\bf y}_j$, which we will do. The statistic can then be used directly to
determine whether to reject $H^j_0$. After substituting ${\bf y}$ by ${\bf
y}_j$,  the approximate MAP statistic (\ref{3-6}) is identical to
\begin{equation} \label{fsv-3}
\frac{  \int (\frac{1}{\sigma})^{T-1}
 e^{-\frac{(T-1) {\hat \sigma}^2_j}{2 \sigma^2} }
d \pi_1(\sigma) }
{\int (\frac{1}{\sigma})^{T-1} e^{-\frac{(T-1) {\hat \sigma}^2_j}{2 \sigma^2}
 -\frac{S_{j}({\hat \phi}_j-\phi_0)^2}{2 \sigma^2}}
d \pi_1(\sigma) }.
\end{equation}

Under the assumption that $\sigma^2$ has a log-normal distribution with mean
$\mu_{\mbox{v}}$ and variance $\tau^2_{\mbox{v}}$ as in Hwang and Liu (2010),
we further approximate (\ref{fsv-3}) by substituting $\sigma$ by its Bayes
estimator; more precisely, $ln(\sigma^2)$
 is substituted by $E(ln(\sigma^2)|\mbox{data})$.
After estimating $\mu_{\mbox{v}}$ and $\tau^2_{\mbox{v}}$ by data in the empirical Bayes fashion,
 we end up with the rejection region
\begin{equation} \label{fsv-4}
 F^j_{\mbox{sv}}=t^2_{jE} > \mbox{crit}, \quad \mbox{where} \quad
  t_{jE}={({\hat \phi}_j -\phi_0)} (\frac{S_{j}}{\hat \sigma_{jE}^2})^{1/2},
\end{equation}
and $\sigma_{jE}^2$ shall be defined below. We note that the critical value,
crit, shall be determined using a bootstrap method so its average Frequentist's
type one error is bounded by $\alpha$. The method is applied to all the other
tests considered in this paper. We do it this way so the proposed tests have
Frequentist's validity.

Note that $t_{jE}$ is the same as $t_{j}$ with the exception that ${\hat \sigma_{j}}^2$ in (\ref{t-test})
is replaced by ${\hat \sigma_{jE}^2}$, which was proposed by Cui et al. (2005).
 To define  ${\hat \sigma_{jE}^2}$,
we take the logarithmic transformation of ${\hat \sigma}^2_j$ and apply the Lindley-James-Stein estimator to estimate
$ln(\sigma^2_j)$. See Lindley (1962) and James and Stein (1961).
We then use the exponential  Lindley-James-Stein estimator
to estimate $\sigma^2_j$. Let $X_j = ln({\hat
    \sigma}^2_j)-E(ln(\chi_{T-2}^2/T-2))$, where
$\chi_{T-2}^2$ is the Chi-squared random variable with degree of freedom $T-2$.
Hence ${\hat \sigma}^2_{jE}$ is the exponential  Lindley-James-Stein estimator,
i.e.
\begin{equation} \label{ts}
{\hat \sigma}^2_{jE}=e^{\delta^j_{LJS}},
\end{equation}
where $$ \delta^j_{LJS}={\bar X}+(1-\frac{(N-3)V_{T-2}}{\sum^N_{j=1}
 (X_j-{\bar X})^2 })_+ (X_j-{\bar X})
 $$
is the Lindley-James-Stein estimator and
  $V_{T-2}$ is the variance of $ln(\chi_{T-2}^2/(T-2))$.

Numerical studies in Section 6 show that (\ref{fsv-4}) has a larger average
power than the t-test. To explain this, note that when ${\hat \sigma}^2_j$ are
very different from each other, $\sum^N_{j=1} (X_j-{\bar X})^2$ is large.
Consequently, $\delta^j_{LJS}$ is close to $X_j$ and hence ${\hat
\sigma}^2_{jE}$ is close to ${\hat \sigma}^2_j$ up to a constant. Therefore
(\ref{fsv-4}) behaves like the t-test and can not be worse. On the other hand,
if ${\hat \sigma}^2_j$ are close to each other, resulting in a small
$\sum^N_{j=1} (X_j-{\bar X})^2$,  $\delta^j_{LJS}$ are close to ${\bar X}$ and
 ${\hat \sigma}^2_{jE}$ are close to the geometric mean of ${\hat \sigma}^2_j$.
 Since ${\sigma}^2_j$ are likely similar to each other, the geometric mean should be a better estimator than
  ${\hat \sigma}^2_j$.
 Hence  test (\ref{fsv-4}) is expected to have a larger average power than the t-test, which is confirmed by
  numerical results.

The Lindley-James-Stein estimator can be derived nonparametrically.
 Hence it is anticipated that the test of
Cui et al. (2005) is robust with respect to the miss-specification of the
distribution of $ln({\sigma}^2_j)$. The conjecture is supported by the
numerical study therein.

\subsection{ The test  shrinking the means: $F_{\mbox{sm}}$}
Now we consider the test that shrinks the means only.
 Case 2 in Section 3 is assumed and hence  Theorem 4 is applicable.

To apply (\ref{3-10}) to the AR(1) model, we consider a normal prior:
\begin{equation}  \label{prior}
\phi_j \sim N(\mbox{$\mu$,  $\tau^2$}) \quad \mbox{when  $H^j_1$ is true.}
\end{equation}
Now we shall evaluate the denominator and the numerator of (\ref{3-10}) where $f(\cdot|\cdot)$ is
defined in (\ref{pdf-1}). The
 denominator of (\ref{3-10}), with ${\bf y}$ being replaced by ${\bf y}_j$, equals
\begin{equation} \label{20-3}
{{ sup_{\sigma}} f({\bf y}_j|\phi_0,\sigma)}
= {sup_{\sigma}} (\frac{1}{\sqrt{2 \pi}\sigma})^{T-1}
e^{-\frac{(T-1) {\hat \sigma}^2_{0j}}{2 \sigma^2}},
\end{equation}
where
\begin{equation} \label{est0}
{\hat \sigma}^2_{0j}=(1/(T-1)) \sum^T_{t=2} (y_{j,t} -\phi_0 y_{j,t-1})^2.
\end{equation}
Direct calculation shows that maximum occurs at $\sigma^2={\hat \sigma}^2_{0j}$. Plugging this into
(\ref{20-3}) shows that (\ref{20-3}) equals
\begin{equation} \label{20-4}
(\frac{1}{\sqrt{2 \pi}{\hat \sigma_{0j}}})^{T-1}  e^{\frac{-(T-1)}{2}}.
\end{equation}
Now to calculate the numerator of (\ref{3-10}), we first calculate
$\int f({\bf y}_j|\phi,\sigma) d \pi(\phi)$
which can be shown after some direct calculations to equal
\begin{equation} \label{20-5}
\left[ (\frac{1}{\sqrt{2\pi}\sigma_j})^{T-1} e^{-\frac{(T-1) {\hat \sigma}^2_j}{2 \sigma^2_j}} \right]
 \cdot
\sqrt{\frac{\sigma^2_j}{S_{j}\tau^2+\sigma^2_j}}
e^{-\frac{1}{2}({\hat \phi}_j-\mu)^2 \frac{S_{j}}{S_{j}
 \tau^2+\sigma^2_j}}.
\end{equation}
The last expression can be derived using (\ref{pdf-2}), rewriting the second exponential term in
 (\ref{pdf-2}) as
 $\frac{\sigma_j}{\sqrt{S_j}}
\frac{\sqrt{S_j}}{\sigma_j} e^{-\frac{1}{2}({\hat \phi}_j-\phi_j)^2
\frac{S_{j}}{\sigma^2_j}} $ and using the classical Bayesian theory which implies that
${\hat \theta}$ has a $N(\mu, \sigma^2+\tau^2)$ distribution if ${\hat \theta}|\theta \sim N(\theta,\sigma^2)$
and  $\theta \sim N(\mu,\tau^2)$.

Now it seems difficult to find the maximum likelihood estimator of ${\sigma}^2_j$. Hence instead we use
${\hat \sigma}^2_j$ (defined in (\ref{sigma})), which maximizes the bracket in (\ref{20-5}).
Hence when ${\hat \sigma}_1$ is taken to be ${\hat \sigma}_j$,  (\ref{3-10}) now equals the ratio
of  (\ref{20-5}) to (\ref{20-4}), which yields
\begin{equation} \label{fsm2}
({\hat \sigma}^2_{0j}/{\hat \sigma}^2_{j})^{(T-1)/2}
\sqrt{\frac{{\hat \sigma}^2_j}{S_{j} \tau^2+{\hat \sigma}^2_j}}
 e^{-\frac{1}{2}({\hat \phi}_j-\mu)^2 \frac{S_{j}}{S_{j} \tau^2+{\hat \sigma}^2_j}}.
\end{equation}

To gain some insight about how  (\ref{fsm2}) works, we show that it can be
approximated by a formula similar to the t-test. Using the approximation
$({\hat \sigma}^2_{0j}/ {\hat \sigma}^2_{j})^{(T-1)/2}= (1+({\hat
\phi}_j-\phi_0)^2 \frac{S_{j}}{(T-1)
 {\hat \sigma}^2_{j}})^{(T-1)/2} \doteq e^{\frac{ S_{j}}{2 {\hat
 \sigma}^2_j} ({\hat \phi}_j-\phi_0)^2}$  when  $T$ is large,
  we express (\ref{fsm2}) as
\begin{equation} \label{fsm3}
\sqrt{\frac{{\hat \sigma}^2_j}{S_{j} \tau^2+{\hat \sigma}^2_j}}
 e^{\frac{ S_{j}}{2 {\hat \sigma}^2_j} ({\hat \phi}_j-\phi_0)^2-\frac{1}{2}({\hat \phi}_j-\mu)^2
 \frac{S_{j}}{S_{j} \tau^2+{\hat \sigma}^2_j}}.
\end{equation}

According to our numerical study, the more important factor in (\ref{fsm3}) is
the exponential part and not the square root factor. Omitting the square root
factor and taking a log transformation of the reminder of (\ref{fsm3}) lead to
\begin{equation} \label{fsm4}
\begin{array}{l}
 {\frac{ S_{j}}{ {\hat \sigma}^2_j} ({\hat \phi}_j-\phi_0)^2-({\hat \phi}_j-\mu)^2
  \frac{S_{j}}{S_{j} \tau^2+{\hat \sigma}^2_j}} \\
 =(\frac{ {\hat \sigma}_j^2 \tau^2} {\tau^2 S_{j} +{\hat \sigma}_j^2})^{-1}
 ( \frac{ \tau^2 S_{j}} {\tau^2 S_{j} +{\hat \sigma}_j^2}
   {\hat \phi}_j + \frac{{\hat \sigma}_j^2 } {\tau^2 S_{j} +{\hat \sigma}_j^2}
 \mu-\phi_0)^2 +   g(S_{j},{\hat \sigma}_j,\mu,\tau) \\
= (\frac { S_{j}} { {\hat \sigma}_j^2 {\hat \beta}_j}  ) ({\hat \phi}^{\star}_j-\phi_0)^2
  +   g(S_{j},{\hat \sigma}_j,\mu,\tau),
  \end{array}
\end{equation}
where
\begin{equation} \label{cort}
{\hat \phi}^{\star}_j =
{\hat \beta}_j  {\hat \phi}_j
+ (1-{\hat \beta}_j) \mu
\quad \mbox{and}  \quad
{\hat \beta}_j= \frac{ \tau^2 S_{j}} {\tau^2 S_{j} +{\hat \sigma}_j^2}.
 \end{equation}
Because  $g(.)$  is a term not involving  ${\hat \phi}_j$, it should be less
relevant to the key parameter $\phi_j$ of the testing problem. Numerical
evidence also suggests that we could ignore the term, which we will do. This
leads to the proposed test, which has a formula similar to the t-test:
\begin{equation} \label{fsm5}
F^j_{\mbox{sm}} = t^2_{jm}, \quad
\mbox{where} \quad t_{jm}=  ({\hat \phi}^{\star}_j-\phi_0)
(\frac { S_{j}} { {\hat \sigma}_j^2 {\hat \beta}_j}  )^{1/2}.
\end{equation}

Note that $F^j_{\mbox{sm}}$ uses the estimator ${\hat \phi}^{\star}_j$ which
shrinks ${\hat \phi}_j$ toward to $\mu$. In application, the hyper-parameters
$\mu$ and $\tau$ are unknown. Hence in Section 5 we use the data to estimate
them in the empirical Bayes fashion.

The formula of  (\ref{fsm5}) works only for $\tau>0$. Later on if $\tau$
is estimated to be zero, (\ref{fsm3}) is used instead.
  This principle applies to all the proposed tests of this paper.

Our numerical studies show that $F_{\mbox{sm}}$  is a reasonable
approximation of (\ref{fsm2}) even when the sample size is as small as
$T=10$. Also, the numerical results in Section 7 indicate that
$F_{\mbox{sm}}$ has higher average power than the t-test.

\subsection{ The test shrinking  the means and the variances: $F_{\mbox{ss}}$}
To produce a test shrinking both means and variances, we assume as in Case 3 of
Section 3 where $(\sigma_j,\phi_j)$ follow the prior distribution
$\pi(\sigma,\phi)=\pi_1(\sigma) \pi_2 (\phi)$, where
  $\pi_2(.)$ is the normal distribution defined in (\ref{prior}), and
$\pi_1(.)$ is the pdf of $\sigma$ with the distribution of $\sigma^2$ being defined
 right after (\ref{fsv-3}).
Applying Theorem 5 and (\ref{map}) to model (\ref{pdf-2}) and
 replacing ${\bf y}$ by ${\bf y}_j$ yields the MAP statistic:
\begin{equation}
\begin{array}{l} \label{fss1}
\frac{ \int \int f({\bf y}_j|\phi,\sigma^2)  d \pi_2(\phi) d \pi_1(\sigma) }
{ \int f({\bf y}_j|\phi_j=\phi_0,\sigma^2) d\pi_1(\sigma)}
\\ =
\frac{ \int (1/\sigma)^{T-1}
 e^{-\frac{(T-1) {\hat \sigma}_j^2}{2 \sigma^2} }
\sqrt{\frac{\sigma^2}{S_{j}\tau^2+\sigma^2}}
 e^{-\frac{1}{2}({\hat \phi}_j-\mu)^2 \frac{S_{j}}{S_{j} \tau^2+\sigma^2}}
d \pi_1(\sigma)
}
{\int (1/\sigma)^{T-1} e^{-\frac{(T-1) {\hat \sigma}_j^2}{2 \sigma^2}-
{\frac{S_{j}({\hat \phi}_j-\phi_0)^2}{2\sigma^2}}
 }
 d \pi_1(\sigma)},
\end{array}
\end{equation}
where the numerator is derived using similar calculations leading to (\ref{20-5}).
Assume  $\pi_1(.)$ is a log normal distribution with mean $\mu_{\mbox{v}}$ and variance $\tau^2_{\mbox{v}}$ as we derived
 $F_{\mbox{sv}}$. Then we can approximate the MAP test by substituting $\sigma^2_j$
 by its Bayes estimator ${\hat \sigma}^2_{jE}$ in the numerator and denominator of (\ref{fss1}), and obtain
 an approximation of the MAP test,
\begin{equation} \label{fss2}
\sqrt{\frac{{{\hat \sigma}_{jE}^2}}{S_{j} {\tau}^2 +{{\hat \sigma}^2_{jE}} }}
e^{  \frac{1}{2}({\hat \phi}_j -\phi_0)^2 \frac{S_{j}} {{\hat \sigma}_{jE}^2}
-\frac{1}{2} ({\hat \phi}_j -\mu)^2  \frac{S_{j}} { S_{j} {\tau}^2+ {{\hat \sigma}_{jE}^2}}}.
\end{equation}

Similar to the calculations leading to  (\ref{fsm4}),
we ignore the first multiple term of (\ref{fss2}) and take a log transformation
of the reminder, yielding
\begin{equation} \label{fss3}
 \begin{array}{c}
 { ({\hat \phi} -\phi_0)^2 \frac{S_{j}} {{\hat \sigma}_{jE}^2}
- ({\hat \phi}_j -\mu)^2  \frac{S_{j}} { S_{j} {\tau}^2+ {{\hat \sigma}_{jE}^2}}} \\
 =
 (\frac { S_{j}} {  {\hat \sigma}_{jE}^2 {\hat \beta}_{jE} } )
 ({\hat \phi}^{\star}_{jE}-\phi_0)^2  + g^{\star}(S_{j},{\hat \sigma}_j,\mu,\tau) ,
\end{array}
\end{equation}
where ${\hat \beta}_{jE}=\frac{ \tau^2 S_{j}} {\tau^2 S_{j} +{\hat \sigma}_{jE}^2}$
and  ${\hat \phi}^{\star}_{jE}={\hat \beta}_{jE}  {\hat \phi}_j + (1-{\hat \beta}_{jE}) \mu$.
Also $g^{\star}(.)$ is a term involving no ${\hat \phi}_j$.
Note we do not need to recalculate (\ref{fss3}) again; we simply replace ${\hat \sigma}_j^2$ by
 ${\hat \sigma}_{jE}^2$ and ${\hat \beta}_{j}$ by ${\hat \beta}_{jE}$ in (\ref{cort}).

Omitting  $g^{\star}(.)$ of (\ref{fss3}) yields the proposed test:
\begin{equation} \label{fss4}
F^j_{\mbox{ss}} = t^2_{jEm}, \quad \mbox{where} \quad
t_{jEm}= ({\hat \phi}^{\star}_{jE}-\phi_0)
 (\frac{ S_{j}} { {\hat \sigma}_{jE}^2 {\hat \beta}_{jE}}   )^{1/2}.
 \end{equation}

The expression of $F^j_{\mbox{ss}}$ not only has a compact formula
similar to the
 t-test, but also enjoys
   nice interpretations. Compared with the t-test,
  $F^j_{\mbox{ss}}$ uses the shrinkage variance estimator ${\hat \sigma}^2_{jE}$
  instead of ${\hat \sigma}^2_{j}$, and  the shrinkage  estimator ${\hat \phi}^{\star}_{j}$
  instead of ${\hat \phi}_{j}$. Therefore,
    $F^j_{\mbox{ss}}$  shrinks the variances as $F_{\mbox{sv}}$ does
     and shrinks the means  as $F_{\mbox{sm}}$ does.
Thus we would expect that $F^j_{\mbox{ss}}$ should perform the best among all the tests.
 Numerical studies in Section 7 confirm this expectation.

Note that we do not need to assume a large $T$ in deriving (\ref{fss2})
 whereas we need it  to derive (\ref{fsm3}). Thus  $F^j_{\mbox{ss}}$ should
 be close to the  MAP test even for small $T$.

\section{The One-sided Test}
We consider the one-sided test: for $1 \leq j \leq N$
\begin{equation}
H^j_0: \quad \phi_j =\phi_0 \quad \mbox{vs.} \quad H^j_1: \quad \phi_j < \phi_0.
\end{equation}
The t-test is to reject $H^j_0$ if $t_j=({\hat \phi}_j -\phi_0) (\frac{S_{j}}
{\hat \sigma_j^2})^{1/2}$ is smaller than a  critical value. To construct tests
having a larger average power, we derive  $RF_{\mbox{sv}}$, $RF_{\mbox{sm}}$
and $RF_{\mbox{sv}}$ which shrink the variances, the means and both the
variances and means, respectively. We include ``R"  in the names of these tests
since $\phi_0$ in the null hypothesis is on the right-hand side of the
alternative region.

\begin{description}
\item{{\it The test shrinking the variances} :  $RF_{\mbox{sv}}$}\\
Suppose  $\phi_j$ is  fixed and unknown, and $\sigma_j$ follows the prior
distribution, defined right after (\ref{fsv-3}), which was used in deriving
$F_{\mbox{sv}}$.  Theorems 1-3 can be directly applied to this problem.
Using  the same arguments leading to $F_{\mbox{sv}}$, we end up with the
test statistic:
\begin{equation}
RF^{j}_{\mbox{sv}}\equiv t_{jE}= ({\hat \phi}_j -\phi_0) (\frac{S_{j}} {\hat \sigma_{jE}^2})^{1/2}.
\end{equation}
The null hypothesis will be rejected if $RF^{j}_{\mbox{sv}}$ is smaller than a critical value.
\item{{\it The test shrinking  the means}: $RF_{\mbox{sm}}$} \\
Suppose $\sigma_j$ is  fixed and unknown and $\phi_j$ is a random variable.
Since the alternative region is $\phi_j<\phi_0$, we postulate that
prior distribution is $N(\mu,\tau^2)$ truncated to the range
 $(-\infty, \phi_0)$.
Hence its pdf  is
\begin{equation} \label{35-1}
f(\phi|\mu,\tau)=\frac{1}{\Phi(\frac{\phi_0-\mu}{\tau})}
 \frac{1}{\sqrt{2\pi}\tau} e^{-\frac{1}{2\tau^2}(\phi-\mu)^2}, \quad \mbox{when} \quad -\infty < \phi < \phi_0,
\end{equation}
where $\Phi(.)$ is the cumulative function of a standard normal distribution.
By Theorem 4, an approximate  MAP test is to reject $H^j_0$ if
\begin{equation} \label{rfsm1}
 \frac{{\mbox{ sup}_{\sigma_j}} \int^{\phi_0}_{-\infty} f({\bf y}_{j}|\phi,\sigma_j^2) d \pi_2(\phi)   }
{{\mbox{ sup}}_{\sigma_j} f({\bf y}_{j}|\phi_j=\phi_0,\sigma_j^2)}  \quad \mbox{is large}.
\end{equation}

Similar to the derivation of $F_{\mbox{sm}}$, a close form for the denominator can be found by replacing
${\hat \sigma}^2_j$ with its maximum point ${\hat \sigma}^2_{0j}$ defined in
(\ref{est0}). However, it does not appear that the numerator has a close form and hence
we simply
replace ${\sigma}^2_{j}$ with ${\hat \sigma}^2_{j}$ defined in
(\ref{sigma}).
 This leads to
\begin{equation} \label{rfsm2}
\begin{array}{c}
 \frac{1}{\Phi(\frac{\phi_0-\mu}{\tau})}
\sqrt{\frac{{\hat \sigma}^2_j}{S_{j} \tau^2+{\hat \sigma}^2_j}}
\Phi (-t_{jm})
  ({\hat \sigma}^2_{0j}/{\hat \sigma}^2_{j})^{(T-1)/2}
 e^{-\frac{1}{2}({\hat \phi}_j-\mu)^2 \frac{S_{j}}{S_{j} \tau^2+{\hat \sigma}^2_j}},
\end{array}
\end{equation}
where  $t_{jm}$ is defined in (\ref{fsm5}).
By adopting the
arguments in  (\ref{fsm3}) and (\ref{fsm4}) for deriving $F_{\mbox{sm}}$,
 the log transformation of $ ({\hat \sigma}^2_{0j}/{\hat
\sigma}^2_{j})^{(T-1)/2}
 e^{\frac{-1}{2}({\hat \phi}_j-\mu)^2 \frac{S_{j}}{S_{j} \tau^2+{\hat
\sigma}^2_j}}$ can be expressed as $(1/2)t_{jm}^2$. Therefore, the MAP
test can be approximated, after  ignoring the first two terms of
(\ref{rfsm2}) and taking the log transformation of the remainder, by
\begin{equation} \label{rfsm3}
log(\Phi(-t_{jm}))+\frac{t^2_{jm}}{2}.
\end{equation}
Using the inequality  $1-\Phi(x)<\frac{1}{x}\phi(x)$ for $x>0$,  (\ref{rfsm3}) can be
shown to be  decreasing in $t_{jm}$ for $t_{jm}>0$. It is obvious that (\ref{rfsm3}) is also decreasing
for $t_{jm}<0$.
Hence (\ref{rfsm3}) is equivalent
to the proposed test which rejects $H^j_0$ when
\begin{equation}
RF^{j}_{\mbox{sm}} \equiv t_{jm}= ({\hat \phi}^{\star}_j-\phi_0)
 (\frac { S_{j}} {{\hat \sigma}_j^2{\hat \beta}_j} )^{1/2} \quad \mbox{is small}.
\end{equation}
\item{{\it The test shrinking the variances and the means:} $RF_{\mbox{ss}}$} \\
Under the assumption of Case 3 in Section 3,
suppose $(\sigma,\phi)$ follows a prior distribution $\pi(\sigma,\phi)=\pi_1(\sigma) \pi_2 (\phi)$,
where $\pi_1(.)$ is the pdf of $\sigma$ with the distribution of
$\sigma^2$ being defined right after  (\ref{fsv-3}) and
 $\pi_2(.)$ is the truncated normal distribution defined in (\ref{35-1}).
Then the MAP test is
\begin{equation} \label{rfss1}
\frac{ \int \int^{\phi_0}_{-\infty} f({\bf y}_j|\phi,\sigma^2)  d \pi_2(\phi) d \pi_1(\sigma) }
{ \int f({\bf y}_j|\phi_j=\phi_0,\sigma^2) d\pi_1(\sigma)}.
\end{equation}
Instead of  integrating with respect to $\sigma$, we replace  $\sigma_j$ in the
 numerator and denominator by ${\hat \sigma}_{jE}$.
An approximation of the MAP test is obtained,
\begin{equation} \label{rfss2}
\begin{array}{c}
 \frac{1}{\Phi(\frac{\phi_0-\mu}{\tau})}
\sqrt{\frac{{\hat \sigma}^2_{jE}}{S_{j} \tau^2+{\hat \sigma}^2_{jE}}}
\Phi (-t_{jEm})
e^{  \frac{1}{2}({\hat \phi} -\phi_0)^2 \frac{S_{j}} {{\hat \sigma}_{jE}^2}
-\frac{1}{2} ({\hat \phi}_j -\mu)^2  \frac{S_{j}} { S_{j} {\tau}^2+ {{\hat \sigma}_{jE}^2}}},
\end{array}
\end{equation}
where   $t_{jEm}$ is defined in (\ref{fss4}).
Ignoring the first two terms of (\ref{rfss2}), taking the log transformation of the remainder, and
using  the leading term in  (\ref{fss3}) to substitute for the exponent yield the proposed statistic:
\begin{equation}
RF^{j}_{\mbox{ss}} \equiv t_{jEm}= ({\hat \phi}^{\star}_j-\phi_0)
(\frac{ S_{j}} {{\hat \sigma}_{jE}^2{\hat \beta}_{jE}})^{1/2}.
\end{equation}

\end{description}

\section{Estimating the Hyper-parameters and the Critical Values}
\subsection{Estimate the hyper-parameters: $\mu$ and $\tau^2$}

\begin{description}
\item{\it The two-sided test: Normal distribution}\\
We follow the empirical Bayes approach and  use data to estimate the  hyper-parameters ($\mu$, $\tau^2$)
used in $F_{\mbox{sm}}$ and  $F_{\mbox{ss}}$.
Suppose
$\phi_j$ follows $N(\mu, \tau^2)$ with probability $\theta_1$, and
 $\phi_j=\phi_0$ with probability $\theta_0=1-\theta_1$. Note that
 $\theta_1$ can be interpreted as $N_1/N$.

By assuming the independence of ${\bf y}_j$ for $j=1,\cdots,N$, the log
likelihood function of ($\mu$, $\tau^2$) is
\begin{equation} \label{like}
\begin{array}{l}
log( f( {\bf y}_1,\cdots, {\bf y}_N|\mu, \tau^2))\\
 = \sum^N_{j=1} log (f({\bf y}_j|\mu, \tau^2)) \\
=
\sum^N_{j=1} log\{\theta_1 \int f({\bf y}_j| \mu, \tau^2, \phi_j=\phi ) d \pi_2 (\phi)  +
(1-\theta_1) f({\bf y}_j|\phi_j=\phi_0)\} \\
 =  C +\sum^N_{j=1} log\{ \theta_1 \sqrt{\frac{\sigma^2_j}{S_{j} \tau^2+\sigma^2_j}}
e^{-\frac{1}{2}  \frac{ S_{j}} { S_{j}\tau^2 +\sigma^2_j}   ({\hat \phi}_j-\mu)^2} +
(1-\theta_1) e^{-\frac{S_j}{2 \sigma^2_j} ( {\hat \phi}_j-\phi_0)^2}\},
\end{array}
\end{equation}
where $C$ is a constant not depending on $\mu$ and $\tau^2$. One can
maximize (\ref{like}) and derive the maximum likelihood estimator for
($\mu$, $\tau^2$, $\theta_1$). However, this involves the maximization
of three variables. Instead,  we propose an  estimator which is easier
to compute. We use the approximation ${\hat \phi}_j |\phi_j, \sigma_j$
$\sim$  $N(\phi_j, \frac{{\sigma}^2_j}{S_{j}})$. Hence
\begin{equation} \label{cond}
\begin{array}{lll}
 E({\hat \phi}_j|\phi_j, \sigma_j ) \doteq \phi_j \quad \mbox{and } \quad
 E({\hat \phi}^2_j|\phi_j, \sigma_j )\doteq \phi^2_j + \frac{{\sigma}^2_j}{S_{j}}.
\end{array}
\end{equation}
Therefore,
\begin{equation} \label{est}
\begin{array}{lll}
 E({\hat \phi}_j) & = & \theta_1 E (\phi_j|\mu,\tau)+(1- \theta_1)\phi_0
 =   \theta_1 \mu + (1-\theta_1) \phi_0 \\
E({\hat \phi}^2_j)& = &  \theta_1  E (\phi^2_j+\frac{\sigma^2_j}{S_{j}} |\mu,\tau)
+(1-\theta_1) (\phi^2_0+E(\frac{\sigma^2_j}{S_{j}})) \\
&= & E(\frac{\sigma^2_j}{S_{j}})+\theta_1(\mu^2+\tau^2)+(1-\theta_1)\phi^2_0.
\end{array}
\end{equation}

Let $m_1$ and $m_2$ denote $E({\hat \phi}_j)$ and $E({\hat \phi}^2_j)$, respectively.
Solving $\mu$ and $\tau^2$ in terms of $m_1$, $m_2$ and $\theta$ yields
\begin{equation} \label{mu}
\mu=\frac{m_1-(1-\theta_1)\phi_0}{\theta_1} \quad \mbox{and} \quad
  \tau^2=\frac{m_2-(1-\theta_1)\phi^2_0 - E(\frac{\sigma^2_j}{S_j})}{\theta_1}-\mu^2.
  \end{equation}
Substitute $m_1$ and $m_2$ with ${\hat m}_1 = (1/N)\sum^N_{j=1} {\hat
 \phi}_j$ and ${\hat m}_2 = (1/N) \sum^N_{j=1} {\hat \phi}^2_j $.
Furthermore, replace  $E(\frac{\sigma^2_j}{S_j})$ in (\ref{mu}) with
$(1/N)\sum^N_{j=1} \sigma^2_j/S_j$, where  $\sigma^2_j$ is, in turn,
replaced with ${\hat \sigma}^2_{jE}$ for  $F_{\mbox{ss}} $ (and ${\hat
\sigma}^2_j$ for  $F_{\mbox{sm}}$). The latter  substitution for
$\sigma^2_j$ is also applied to (\ref{like}). Moreover, plug the resultant
formula for $\mu$ and $\tau^2$ into (\ref{like}). Then the resultant pseudo
likelihood function is a function of $\theta_1$ only. We then estimate
$\theta_1$ by  ${\hat \theta}_1$ which maximizes the function. Using ${\hat
\theta}_1$,  ${\hat m}_1$ and ${\hat m}_2$, we may estimate $\mu$  and
$\tau^2$ based on (\ref{mu}).

\item{\it The one-sided test: Truncated normal distribution}\\
To estimate the hyper-parameter ($\mu$,$\tau^2$) of a truncated normal distribution used to derive
 $RF_{\mbox{sm}}$ or $RF_{\mbox{ss}}$, we adopt the empirical Bayes
approach again, assuming that $\phi_j$ follows the truncated  normal
distribution with probability $\theta_1$ and  $\phi_j=\phi_0$ with
probability $1-\theta_1$. Therefore, the log  likelihood function of
($\mu$, $\tau^2$) is
\begin{equation} \label{lik2}
\begin{array}{l}
log( f( {\bf y}_1,\cdots, {\bf y}_N|\mu, \tau^2) ) \\
 =  C +\sum^N_{j=1} log\{ \theta_1 \frac{1}{\Phi(\frac{\phi_0 -\mu}{\tau})}
 \Phi(-t^o_{jm})
\sqrt{\frac{\sigma^2_j}{S_{j} \tau^2 +\sigma^2_j}}
e^{-\frac{1}{2}  \frac{ S_{j}} { S_{j}\tau^2 +\sigma^2_j}   ({\hat \phi}_j-\mu)^2} +
(1-\theta_1) e^{-\frac{S_j}{2 \sigma^2_j} ( {\hat \phi}_j-\phi_0)^2}\},
\end{array}
\end{equation}
where $t^o_{jm}$ is identical to  $t_{jm}$ in but using $\sigma_j$ to replace ${\hat \sigma}_j$.

Using (\ref{cond}) and the moments of a truncated normal distribution, we
have the following results after some calculations:
\begin{equation}
\begin{array}{lll} \label{trmu}
E({\hat \phi}_j)  & = & \theta_1(\mu- \lambda(\alpha) \tau)+(1-\theta_1) \phi_0 \\
E({\hat \phi}^2_j) & = &  E(\frac{{ \sigma}^2_j}{S_{j}}) +
                   \theta_1 \{\tau^2(1-\alpha \delta(\alpha))+(\mu-\lambda(\alpha)\tau)^2 \} +
                   (1-\theta_1) \phi^2_0
\end{array}
\end{equation}
where $\alpha=\frac{\phi_0-\mu}{\tau}$, $\lambda(\alpha)=\frac{\phi(\alpha)}{\Phi(\alpha)}$ and
$\delta(\alpha)=\lambda(\alpha)(\alpha+\lambda(\alpha))$.
Replacing $E({\hat \phi}_j)$ and
 $E({\hat \phi}^2_j)$ by
 $m_1$ and $m_2$ in (\ref{trmu}), respectively, gives us
\begin{equation} \label{mu2}
\begin{array}{l}
\mu  =  \frac{m_1 -(1-\theta_1)\phi_0}{\theta_1}+\lambda(\alpha)\tau \\
\tau^2= \{\frac{m_2-(1-\theta_1)\phi^2_0  - E(\frac{\sigma^2_j}{S_j})}{\theta_1} -
(\mu-\lambda(\alpha)\tau)^2\}/(1-\alpha\delta(\alpha)).
\end{array}
\end{equation}
Since the right-hand side of  (\ref{mu2}) still involves $\mu$ and
$\tau^2$, an iteration algorithm is proposed to estimate ($\mu$,$\tau^2$).
We use the estimator for ($\mu$, $\tau^2$) in the two-sided case depicted
above as the initial value to obtain a function of $\theta_1$ only.
Calculate ${\hat \theta}_1$ that maximizes the function. Now plug  ${\hat
\theta}_1$  and the initial value of ($\mu$, $\tau^2$) into the right-hand
side of (\ref{mu2}) to obtain a new estimator of $\mu$ and $\tau^2$. The
process is repeated to obtain a new estimator of $\theta_1$, $\mu$ and
$\tau^2$. In the above calculation $m_1$ and $m_2$ are replaced by
 $(1/N)\sum^N_{j=1} {\hat \phi}_j$ and $(1/N) \sum^N_{j=1} {\hat
 \phi}^2_j$, respectively. Also  $E(\frac{\sigma^2_j}{S_j})$ is
replaced by $(1/N)\sum^N_{j=1} \sigma^2_j/S_j$, where $\sigma^2_j$ is,
in turn, replaced by ${\hat \sigma}^2_{jE}$ for  $F_{\mbox{ss}}$ (and
${\hat \sigma}^2_{j}$ for $F_{\mbox{sm}}$). The latter substitutions
for $\sigma_j^2$ are also applied to (\ref{lik2}).
\end{description}

\subsection{Generating the critical value by the Bootstrap method}
In order to have a good finite sample property, we should use the bootstrap method
to determine the critical values of the proposed tests.
In what follows, we present the details of the bootstrap procedure for the two-sided test.
A similar procedure can be applied to the one-sided test.

Let
\begin{equation} \label{boot-1}
{\hat e}_{j,t}=y_{j,t}- {\hat \phi}_j y_{j,t-1} \quad \mbox{for $2 \leq t \leq T$ and
  $1 \leq j \leq N$.}
\end{equation}
Under the null hypothesis, we use the hypothesized value $\phi_0$ to create the
following bootstrap sample for the $j$-th group,
$\{y^{\star}_{j,t},t=1,2,,\cdots,T \}$, where
\begin{equation}  \label{boot-2}
 y^{\star}_{j,t} = \phi_0 y^{\star}_{j,t-1} +e^{\star}_{j,t}, \quad t=1,2, \cdots, T,
 \end{equation}
 and  ${e^{\star}_{j,t}}'s$  are sampled  with replacement from
 $\{{\hat e}_{j,t}, \mbox{$2 \leq t \leq T$}\}$.

One can plug  ${y^{\star}_{j,t}}'s$ into the  t-statistic,  $F_{\mbox{sv}}$ statistic,
 $F_{\mbox{sm}}$ statistic and $F_{\mbox{ss}}$ statistic.
For each statistic, repeat it $R$ times and calculate the percentile ($95\%$tile for $5\%$ test)
which is then used as the critical value.

Note that in calculating the critical values for  $F_{\mbox{sm}}$ and
$F_{\mbox{ss}}$, we use data to estimate $\mu$ and $\tau$ once and from then
on, $\mu$ and $\tau$ are set to be identical to its estimated value.  Hence in
each bootstrap sample, $\mu$ and $\tau$ are not re-estimated. This is
reasonable since in the bootstrap sample, $\phi_j$ is taken to be $\phi_0$, the
hypothesized value. The bootstrap samples do not have information about
$\phi_j$ and hence they should not be used to estimate the hyper-parameters of
$\phi_j$. Regarding ${\hat \sigma}^2_{jE}$ used in the two tests
$F_{\mbox{sv}}$ and $F_{\mbox{ss}}$, we do recalculate its value for each
bootstrap sample, since they contain the information about $\sigma^2_j$.

\section{Simulation Studies}
\subsection{ The white noise hypothesis: Two-sided test for $\phi_0=0$ }
This simulation considers a special case of the two-sided test in which the
null and the alternative hypotheses are, respectively, $H^j_0$: $\phi_j=0$ and
$H^j_1$:  $\phi_j \neq 0$
 for $1 \leq j \leq N$. The null hypothesis is  commonly referred to as the white noise hypothesis.

In Section 6, we estimate the hyper-parameters of the proposed tests under the
assumption  that cross section series are mutually independent; namely, the $N$
series are independent. However,  the following simulation studies both the
independent and dependent cases. In general, we assume a multi-factor error
structure (Pesaran et al., 2013), which includes both independent and dependent
cases, in order to check the robustness of the proposed tests with respect to
cross section dependence.
  It turns out that in both the independent case and dependent case,
  the proposed tests improve uniformly  over the t-test. Hence in both cases, the proposed tests apparently
 ``borrow the strength" from all the populations to do better than the t-test.

Specifically,  the data are generated using the model
\begin{equation} \label{model}
 \mbox{$y_{j,t}=\phi_j y_{j,t-1}+e_{j,t}$ \quad for $ 1 \leq j \leq N$,
 with $\phi_j=0$ for $j >N_1$},
\end{equation} \label{err}
\begin{equation} \label{factor}
e_{j,t}=c_{j,1}f_{1,t}+c_{j,2}f_{2,t}+\epsilon_{j,t},
\end{equation}
where $\epsilon_{j,t}$, $1 \leq j \leq N$, $1 \leq t \leq T$, are independently
 $N(0,\sigma^2_j)$ distributed.

 Model (\ref{factor}) is called a
multi-factor model, which reduces to the independent model when
$c_{j,1}=c_{j,2}=0$. Otherwise, $\{ e_{j,t} \}$, for $1 \leq j \leq N$, are
dependent. For the dependent case studied below, $c_{j,1}$ and $c_{j,2}$ are
generated as random samples from the uniform distribution over $[0,1]$ and
$[0,2]$ respectively.

For t-statistic, $F_{\mbox{sv}}$, $F_{\mbox{sm}}$ and $F_{\mbox{ss}}$, we then
calculate the average power $\mbox{ETP}/N_1$  and the average type one error
$\mbox{EFP}/(N-N_1)$, where $\mbox{ETP}$ and $\mbox{EFP}$
 are defined in
(\ref{etp-1}) and (\ref{efp-1}).

 The parameters $\sigma_j$, $1 \leq j \leq N$,
 are i.i.d samples generated from $\pi_1(\sigma)$ and $\phi_j$, $1 \leq j \leq N_1$,
 are i.i.d samples generated from $\pi_2(\phi)$, where $\pi_1$ and $\pi_2$ will be specified below.

\begin{description}
\item{\it  Normal prior distributions} \\
The prior distributions assumed are
\begin{equation} \label{prior-mt}
\begin{array}{l}
\pi_1(\sigma): \quad ln(\sigma^2_j) \sim N(\mu_{\mbox{v}}, \tau^2_{\mbox{v}}) \quad \mbox{for}
 \quad 1 \leq j \leq N, \\
\pi_2(\phi): \quad \phi_j \sim N(\mu, \tau^2) \quad \mbox{for}  \quad 1 \leq j \leq N_1.
\end{array}
\end{equation}

We now examine the average power and the average type one error of the
t-test and our proposed tests. In each of Figures 1.1 through 1.6, the
simulated average power and the  average type one error are plotted,
against $\mu$, in  solid curve and dotted curve respectively. In the
simulation, each  point is based on averaging at least 4,000
 replications.

 In Figures 1.1 through 1.5, the cross section series are mutually
 independent. For various settings of $T$, $N$ and $N_1$, $\tau$ and the
 coefficient of variation ($\mbox{CV}=\tau_{\mbox{v}}/\mu_{\mbox{v}}$)
specified in the headings of these figures, the figures basically show that
\begin{equation} \label{statment}
\begin{array}{l}
\mbox{{\bf statement} (i)}:
\begin{array}{l}
\mbox{all the proposed tests have uniformly higher average  power than} \\ \mbox{the
t-test,}        \end{array}
 \\
\mbox{{\bf statement} (ii)}: \mbox{the uniformly most powerful test is $F_{\mbox{ss}}$ test.}
\end{array}
\end{equation}
 And, the average power of $F_{\mbox{ss}}$ could be $70 \%$, as shown in
Figures 1.1,  larger than that of the t-test.

Moreover all the tests have average type one error controlled under $5 \%$
level with the exception of Figure 1.5,  which correspond to small $N$ and
$N_1$. Further numerical study shows that the discrepancy is due to the
estimation error of $\mu$ and $\tau$, which is larger for small $N$ and
$N_1$. However, even in Figure 1.5, the average type one errors  of
alternative tests are only slightly larger than 0.05.

As for the case of cross section dependence, we adopt the multi-factor
model (\ref{factor}) to generate the data. Under the same settings of $T$,
$N$ and $N_1$, $\tau$ and $\mbox{CV}$ as those in Figures 1.1 through 1.5,
we obtain very similar graphs showing basically that the statements (i) and
(ii) in (\ref{statment}) hold. We only report Figure 1.6 having the setting
of Figure 1.1.

In fact, the study shows that the improvements of $F_{\mbox{ss}}$ test over
the t-test are slightly larger in some of the dependent cases. This is
intuitively reasonable since a procedure shrinking toward the  common means
or variances should be expected to do better when the sections are more
correlated.

Our simulation studies also confirm the effectiveness of the estimator for
the hyper-parameter ($\mu,\tau^2)$ in Section 6.1. More specifically, the
average power of the proposed tests using the estimated $(\mu, \tau^2)$ is
very similar to that of the tests using the true values, although the
average power corresponding to the true values is not reported here.

In Econometrics, it is important to focus on the alternative hypothesis
which is close to the null hypothesis. This is especially true for the unit
root test, to be discussed in Section 7.2. Consequently, the tests do not
have large average  power. However,  the increase of the average power by
$0.05$ will, on the average, increase the detected true positives by
$(0.05)N$, which could be quite substantial when $N$ is large.

\item{\it  Uniform prior distributions and fixed effect model} \\
To show the robustness of the proposed tests with respect to the
 miss-specification of prior distributions, we use ``wrong" distributions
such as the  uniform distributions and fixed effect model to generate
parameters. We consider the  uniform distributions as
\begin{equation}
\begin{array}{l} \label{unfm}
\pi_1(\sigma): \quad ln(\sigma^2_j) \sim U(2-2\sqrt{3}\tau_{\mbox{v}}, 2+2\sqrt{3}
 \tau_{\mbox{v}} ) \quad \mbox{for}
\quad 1 \leq j \leq N, \\
\pi_2(\phi): \quad \phi_j \sim U(\mu-2\tau, \mu+2\tau) \quad \mbox{for}
 \quad 1 \leq j \leq N_1.
\end{array}
\end{equation}
We write the distribution of $ln (\sigma^2_j)$ this way, so that
the variance is $4\tau^2_{\mbox{v}}$ and the mean is two; consequently
 $\mbox{CV}=\tau_{\mbox{v}}$. For such a prior, we plot the average power
 using the same settings as Figures 1.1 through 1.5 for both independent
 case and multi-factor models. The resultant graphs are similar to Figures
 1.1 through 1.5. We report only Figure 2.1 (corresponding to the
 independent case) and Figure 2.2 (corresponding to the multi-factor model)
 both having the same settings as Figures 1.1. Both figures and the
 unreported figures  basically confirm the two statements in
 (\ref{statment}).

To study the fixed effect model,  i.e. $\phi_j$ being fixed, let
\begin{equation} \label{fix}
\begin{array}{l}
 \phi_j =\mu-2\tau \quad \mbox{for} \quad 1 \leq j \leq {N_1}/{2},  \\
 \phi_j =\mu + 2\tau  \quad \mbox{for} \quad {N_1}/{2}+1 \leq j \leq N_1,
 \end{array}
\end{equation}
and $\sigma_j=\sigma$ for all $j$ (CV=0). The results displayed in Figures
2.3 and 2.4 show that the improvements obtained by the proposed tests are
also robust with respect to this ``wrong"  setting. Statements in
(\ref{statment}) basically hold.

\item{\it  Conditional heteroscedasticity} \\
Below, we shall  generate $\epsilon_{j,t}$ from a GARCH(1,1) model instead
of an i.i.d Normal model. The GARCH(1,1) model is commonly used in Finance
and Economics to  describe conditionally heteroscedastic phenomena. The
GARCH(1,1) model used is
\begin{equation} \label{garch}
\epsilon_{j,t}=\omega_{j,t} \epsilon^{\star}_{j,t} \quad \mbox{where}  \quad
{\omega^2_{j,t}}=1+0.8 {\omega^2_{j,t-1}}  + 0.15 \epsilon^2_{j,t-1};
\end{equation}
where $\epsilon^{\star}_{j,t}$ are i.i.d. standard normal random variables
and $w^2_{j,t}$ is the conditional variance of $\epsilon_{j,t}$. Models
(\ref{model}), (\ref{factor}) and (\ref{prior-mt}) are still assumed except
$\epsilon_{j,t}$ follows (\ref{garch}). The results in Figures 3.1 and 3.2
show that the proposed tests still improve on the t-test. In particular,
Figure 3.2 assumes a model that has cross section dependence and
conditional heteroscedasticity. Therefore the improvements are quite robust
with respect to dependence and heteroscedasticity. Statements in
(\ref{statment}) are basically correct.

\item{\it  Large dimensional series} \\
In what follows, we study the large dimensional series that $N=1000$ and
$N_1=500$. Under the settings of Figure 1.1, Figures 4.1 and 4.2 report the
results of independent and dependent cases, respectively. Both figures
strongly confirm the two statements in (\ref{statment}) and show that the
improvements provided by the proposed tests, including $F_{\mbox{sv}}$,
$F_{\mbox{sm}}$ and $F_{\mbox{ss}}$, over the t-test increase slightly when
the dimension increases.

\end{description}

\subsection{The unit root hypothesis: One-sided test for $\phi_0=1$}
Now we apply all the tests to the unit root hypothesis, for testing  $H^j_0$:
$\phi_j=1$ vs. $H^j_1$:  $\phi_j < 1$ for $1 \leq j \leq N$.

We generate the data using the model
\begin{equation}
 \mbox{$y_{j,t}=\phi_j y_{j,t-1}+e_{j,t}$ \quad for $ 1 \leq j \leq N$,
 with $\phi_j=1$ for $j >N_1$},
\end{equation}
where  ${e_{j,t}}'s$ are generated from equation (\ref{factor}),
${\epsilon_{j,t}}'s$ from a normal distribution $N(0,\sigma^2_j)$,
${\sigma^2_{j}}'s$ from a  prior distribution $\pi_1(\sigma)$ for $1 \leq j
\leq N, $ and $\phi_j$ from $\pi_2(\phi)$  for $1 \leq j \leq N_1$, where
$\pi_1$ and $\pi_2$ are specified below.

 We shall graph the  average powers and the average type one errors
 of  t-test, $RF_{\mbox{sv}}$ and  $RF_{\mbox{ss}}$.
However, we do not show the results of  $RF_{\mbox{sm}}$ since its performance
is very similar to (but slightly worse than)  the t-test.

\begin{description}
\item{\it  Truncated normal prior distributions} \\
We generate parameters by the prior distributions  to derive the proposed
tests. That is, $ln(\sigma^2_j)$ for $1 \leq j \leq N$ have
$N(\mu_{\mbox{v}}, \tau^2_{\mbox{v}})$ distribution, and $\phi_j$, $1 \leq
j \leq N_1$, follow  a $N(\mu, \tau^2)$ distribution truncated to the range
 $(-\infty,1)$.

In Figures 5.1 through 5.5, we graph the simulated average power (plotted
by solid lines), and the simulated average type one error (plotted by
dotted lines) of the three tests under the various combinations of $T$,
$N$, $N_1$, CV and $\tau$ specified in the headings. These figures deal
with the cases when the cross sections are independent, namely model
(\ref{factor}) with $c_{j,1}=c_{j,2}=0$ for all $j$. These graphs
demonstrate that the following statement (iii) holds:
\begin{equation} \label{stat2}
\begin{array}{l}
\mbox{{\bf statement} (iii)}:
\begin{array}{l} \mbox{$RF_{\mbox{ss}}$ and  $RF_{\mbox{sv}}$ basially have uniformly
 higher average power} \\ \mbox{ than the t-test.} \end{array}
\end{array}
\end{equation}
Hence statement (i) in (\ref{statment}) basically holds with
$RF_{\mbox{sv}}$ and $RF_{\mbox{ss}}$. Regarding the question as to which
test of the two is better, the answer is not clear. In principle, the test
$RF_{\mbox{ss}}$ should perform better since it has more to do with the
specifics of the priors. However, $RF_{\mbox{ss}}$ is not  always the
winner. This may have to do with the fact that more hyper-parameters need
to be estimated in constructing $RF_{\mbox{ss}}$ than  those in
constructing $RF_{\mbox{sv}}$. Estimation of the hyper-parameters is  a
difficult problem in the one-sided case because of truncation of the prior.
This may explain why $RF_{\mbox{ss}}$ is not always the winner.

 For the dependent case, we also produce results similar to Figures 5.1 to
5.5. However, only Figure 5.6 is reported which has the same settings as in
Figure 5.1. Figure 5.6, for the dependent multi-factor model, is very
similar to Figure 5.1 for the independent model. This demonstrates that the
 improvements of the proposed tests over
the t-test are quite robust with respect to the cross section dependence.

Whether one uses $RF_{\mbox{sv}}$ or $RF_{\mbox{ss}}$,  these figures show
that both  have higher average power than the t-test. The average power of
$RF_{\mbox{ss}}$ could be about $25 \%$ larger than the t-test (Figures 5.1
and 5.6). The average type one error of all the tests are  controlled under
or nearly under the $5\%$ level.

\item{\it  Uniform prior distributions and fixed effect model}\\
Below we shall study different priors and models. In all cases, statement
(iii) in (\ref{stat2}) is shown to be true. Specifically, $RF_{\mbox{sv}}$
and $RF_{\mbox{ss}}$ have uniformly greater power than t-test and there is
no clear winner between $RF_{\mbox{sv}}$ and $RF_{\mbox{ss}}$.

To study how improvements are affected by a ''wrong prior", we consider
 (\ref{unfm})  except that ${\phi_j}$ is truncated so  ${\phi_j}$ is in
$(-\infty, 1)$ for $1 \leq j \leq N_1$. Following the settings of Figures
5.1 through 5.5, we plot the average powers which show that the fact that
the prior is the ''wrong" prior has little effect. The resultant figures
are very similar. We only report Figure 6.1 (similar to Figure 5.1) and
Figure 6.2 (similar to Figure 5.6).

Similar plottings were carried out for a fixed effect model where $\phi_j
=\mu-2\tau$ for $1 \leq j \leq {N_1}/{2}$ and $\phi_j =min(0.99,\mu +
2\tau)$  for ${N_1}/{2}+1 \leq j \leq N_1$. Since in the graphs, $\mu \leq
1$, the choice of $\phi_j$ above ensures that $\phi_j<1$, for $1 \leq j
\leq N_1$. Hence the first $N_1$ hypotheses are the alternative hypotheses.
The rest of the hypotheses are the null hypotheses, where $\phi_j=1$ for $j
>N_1$.

We report the results in Figures 6.3 and 6.4 which have the same settings
as Figures 6.1 and 6.2, respectively. These two sets of graphs are very
similar, confirming statement (iii).

\item{\it  Conditional heteroscedasticity} \\
 Figures 7.1 and 7.2 graph the average powers and average type one errors
when the data are  generated by equation (64) with the GARCH(1,1) error
(\ref{garch}), and the parameters are generated by the truncated normal
prior distributions. The results confirm statement (iii) and show that the
proposed tests still improve on the t-test even when conditional
heteroscedasticity and cross section dependence are present.

\item{\it  Large dimensional series} \\
Figures 8.1 and 8.2 demonstrate results when the dimensions $N=800$ and
$N_1=600$. The figures show that statement (iii) still holds when dimension
is large.
\end{description}

\section{Concluding Remarks}
To analyze the coefficients of a panel AR(1) model, we propose tests which
 determine which individual hypothesis should be accepted or rejected.
Furthermore,  our proposed tests improve on the  t-test under the criterion of
average power. We derive them using empirical Bayes approach and then using
approximation to obtain our proposed tests, which have a form similar to the
t-test. The only difference is that, in our proposed tests, the estimators of
the means and variances are replaced by shrinkage estimators. The proposed
tests ``borrow the strength" from all the series to test against every
individual series, resulting in more power. Simulation studies show that the
proposed tests have significant improvements over the t-test, especially
  when the sample size $T$ is small and the dimension $N$ is moderate or large.
  Compared to the t-test, the average power of $F_{\mbox{ss}}$ and
$RF_{\mbox{ss}}$ could be  $70\%$ higher in the two-sided test, and $25\%$
higher in the one-sided test respectively.

In this paper, we derive the tests under the assumption that the series are
independent; and show that  ``borrowing the strength'' from independent series
will improve average power of the t-test. However, simulation demonstrates that
the improvement is robust with respect to the cross section dependence. This is
only reasonable. A procedure that can do well by ``borrowing the strength" even
from independent series can certainly do so from dependent series. In this
paper, we only work with AR(1) model; we, however, anticipate that these
results can be generalized to the other more complex time series models.
 Since the proposed tests can determine acceptance or rejection of an
individual hypothesis, this should prove to be a  very useful method in
practice.

${}$\\
\begin{description}
\item{\Large{\bf Appendix}}
\end{description}
\begin{description}
\item{\bf Proof of Theorem 1:} The difference in (\ref{3-1}) equals
\begin{equation} \label{a-1}
\frac{1}{N_1} \sum_{\{j|\phi_j \in D\}}[ g_j(\sigma_j) -E( g_j(\sigma_j))],
\end{equation}
where  $g_j(\sigma_j) = P_{\phi_j, \sigma_j}({\bf y}_j \in C) = \int_{{\bf
 y} \in C} f({\bf y}| {\phi_j, \sigma_j}) d{\bf y}$ and $E(g_j(\sigma_j))=
 \int g_j(\sigma)d \pi_1(\sigma)$. Since variance of $g_j(\sigma_j)$ $\leq$
 $E(g^2_j(\sigma_j)) \leq 1$, the variance of (\ref{a-1}) is bounded above
 by $p_1/p^2_1=1/p_1$, which converges to zero. Hence (\ref{a-1}) converges
 in probability to zero by the law of large numbers, completing the proof
 of (\ref{3-1}).

 Equation (\ref{3-2}) can be proved similar, except noting that $\int
 \int_{{\bf y} \in C} f({\bf y}| {\phi_0, \sigma_j}) d{\bf y} d \pi_1
 (\sigma_j)$
does not depend on $j$.
\item{\bf Proof of Theorem 3:}
It sufficient to show that
\begin{equation}\label{a-2}
\mbox{sup}_{\phi_1\cdots,\phi_N \in D} {\frac{1}{N_1} \sum_{\{j|\phi_j \in D\}}
\int  f({\bf y}| {\phi_j, \sigma}) d\pi_1 (\sigma)} =
\int  f({\bf y}| {{\hat \phi}_M, \sigma}) d\pi_1 (\sigma).
\end{equation}
Note that the left-hand side is obviously bounded by the right-hand side since
 $f({\bf y}| {\phi_j, \sigma}) \leq  f({\bf y}| {{\hat \phi}_M, \sigma})$ for $\phi_j \in D$.
 Also replacing $\phi_j$ by ${\hat \phi}_M$ on the left-hand side leads to a lower bound, which is
 exactly the right-hand side, establishing (\ref{a-2}) and the theorem.
\end{description}

${}$\\
\begin{description}
\item{\Large{\bf  Acknowledgements}}
\end{description}
The research of  J. T. Gene Hwang is supported by the Ministry of Science and
Technology, Taiwan in grant NSC-100-2118-M-194-004-MY3, as well as by a
research grant from the Foundation for the Advancement of Outstanding
Scholarship, Taiwan. The research  of Yu-Pin Hu  is supported by the Ministry
of Science and Technology, Taiwan in grant NSC-103-2118-M-260-001. The authors
also want to thank Ms. Rebecca Brody and Ms. Jinghuei Hwang for their careful
editing of the paper.

${}$\\
\begin{description}
\item{\Large{\bf References}}
\end{description}
\begin{description}
\item Baltagi, B.H.,  Kao, C., 2000.
Nonstationary panels, cointegration in panels and dynamic panels. A survey.
 In: Batltagi, B. (Ed.) Nonstationary Panels, Panel Cointegration, and
 Dynamic Panels. In: Advances in Econometrics, vol. 15. JAI Press,
 Amsterdam, pp. 7-52.

\item Bai, J., Ng, S., 2004.
A panic attack on unit root tests and cointegration. Econometrica  72,
1127-1177.

\item Bai, J.,  Ng, S., 2010,
Panel unit   root tests with cross-section dependence: A further
 investigation.  Econometric Theory  26,  1088-1114.


\item Cui, X., Hwang, J., Qiu, J., Blades, N.J.,  Churchill, A., 2005.
Improved statistical tests for differential gene expression by shrinking
variance components estimates.  Biostatistics  6, 59-75.

\item Ferretti, N., Romo, J., 1996. Unit root bootstrap tests for AR(1) models.
 Biometrika 83, 849-860.

\item Hwang, J., Liu, P., 2010.
Optimal testing shrinking both means and variances applicable to microarray
data.  Statistical Applications in Genetics and Molecular Biology  9,
Article 36.

\item James, W., Stein, C., 1961.
Estimation with quadratic loss.  Proceeding Fourth Berkeley Symposium  on
Mathematical Statistics and Probability 1, 361-379. Berkeley: University of
California Press.

\item Levin, A.,  Lin, C., Chu, C.J.,  2002.
Unit root tests in panel data: asymptotic and finite sample properties.
Journal of Econometrics 108, 1-24.

\item Lindley, D.V., 1962.
Discussion of Professor Stein's paper: Confidence sets for the mean of a
multivariate normal distribution.  Journal of the Royal Statistical Society
B  24, 285-287.

\item Liu, P., 2006.
Sample size calculation and empirical Bayes test for microarray data. Ph.D.
Dissertation, Cornell University.

\item Noma, H., Matsui, S., 2012.
The optimal discovery procedure in multiple significance testing: an
empirical Bayes approach.  Statistics in Medicine 31,  165-176.

\item Noma, H., Matsui, S., 2013.
An empirical Bayes optimal discovery procedure based on semiparametric
hierarchical mixture models. Computational and Mathematical method in
Medicine, Article 568480.

\item Park, J.Y., 2003.
Bootstrap unit root tests.  Econometrica  71,  1845-1895.

\item Pesaran, M.H., 2007.
A simple panel unit root test in the presence of cross section dependence.
Journal of Applied Econometrics 22, 265-312.

\item Pesaran, M.H., Smith, L.V., Yamagata, T., 2013.
Panel unit root tests in the presence of a multifactor error structure.
Journal of Econometrics 175, 94-115.

\item Philips, C.B., Xiao, Z., 1998.
A primer on unit root testing. Journal of Economic Surveys 12, 423-469.

\item Smyth, G.K., 2004. Linear models and empirical Bayes methods for assessing differential expression
in microarray experiments. Statistical Applications in Genetics and
Molecular Biology 3, Article 3.

\item Storey, J.D., 2007.
The optimal discovery procedure: A new approach to simultaneous
significance testing.   Journal of the Royal Statistical Society B  69,
347-368.

\item Storey J.D., Dai, J.Y.,  Leek, J.T., 2007.
The optimal discovery procedure for large-scale significance testing, with
applications to comparative microarray experiments.  Biostatistics  8,
414-432.

\item Uhlig, H., 1994.
What macroeconomists should know about unit roots:  A Bayesian perspective.
 Econometrics Theory 10, 645-671.
\end{description}

\end{document}